\documentclass{article}[12pt]
\usepackage{amssymb,amsfonts}
\usepackage[pdftex]{graphicx}
\usepackage{here}

\begin{document}

\begin{center}
{\Large\bf The Maillet--Malgrange type theorem for generalized power series}
\end{center}

\begin{center}
\large \bf R.\,R.\,Gontsov, I.\,V.\,Goryuchkina
\end{center}
\bigskip

\begin{abstract}
There is proposed the Maillet--Malgrange type theorem for a generalized power series (having complex power exponents) 
formally satisfying an algebraic ordinary differential equation. The theorem describes the growth of the series coefficients.  
\end{abstract}

\section{Introduction} 

Let us consider an ordinary differential equation (ODE) 
\begin{eqnarray}\label{ADE}
F(z,u,\delta u,\ldots,\delta^m u)=0
\end{eqnarray}
of order $m$ with respect to the unknown $u$, where $F(z,u_0,u_1,\ldots,u_m)\not\equiv0$ is a polynomial of $m+2$ variables, $\delta=z\frac d{dz}$.

The classical Maillet theorem \cite{Mai} asserts that any formal power series solution
$\varphi=\sum_{n=0}^{\infty}c_nz^n\in{\mathbb C}[[z]]$ of (\ref{ADE}) 
is a power series of Gevrey order $1/k$ for some $k\in{\mathbb R}_{>0}\cup\{\infty\}$. This means that 
the power series
$$
f=\sum_{n=0}^{\infty}\frac{c_n}{\Gamma(1+n/k)}\,z^n
$$
converges in a neighbourhood of zero, where $\Gamma$ is the Euler gamma-function. In other words,
$$
|c_n|\leqslant A\,B^n\,(n!)^{1/k}
$$
for some $A,B>0$. 

First exact estimates for the Gevrey order of a power series formally satisfying an ODE were obtained by J.-P.\,Ramis \cite{Ram} in the linear case,
$$
Lu=a_m(z)\delta^mu+a_{m-1}(z)\delta^{m-1}u+\ldots+a_0(z)u=0, \qquad a_i\in{\mathbb C}\{z\}.
$$
He has proved that such a power series is of the {\it exact} Gevrey order $1/k\in\{0,1/k_1,\ldots,1/k_s\}$, where $k_1<\ldots<k_s<\infty$ are all of the positive slopes of the Newton polygon ${\cal N}(L)$ of the operator $L$. 
This polygon is defined as the boundary curve of a convex hull of a union $\bigcup_{i=0}^mX_i$, 
$$
X_i=\{(x,y)\in{\mathbb R}^2\mid x\leqslant i,\;y\geqslant{\rm ord}_0a_i(z)\}, \qquad i=0,1,\ldots,m
$$
(see Fig. 1). The exactness of the Gevrey order means that there is no $k'>k$ such that the power series is of the Gevrey order $1/k'$.

\begin{figure}[H]
 \centering
\includegraphics*[height=4cm]{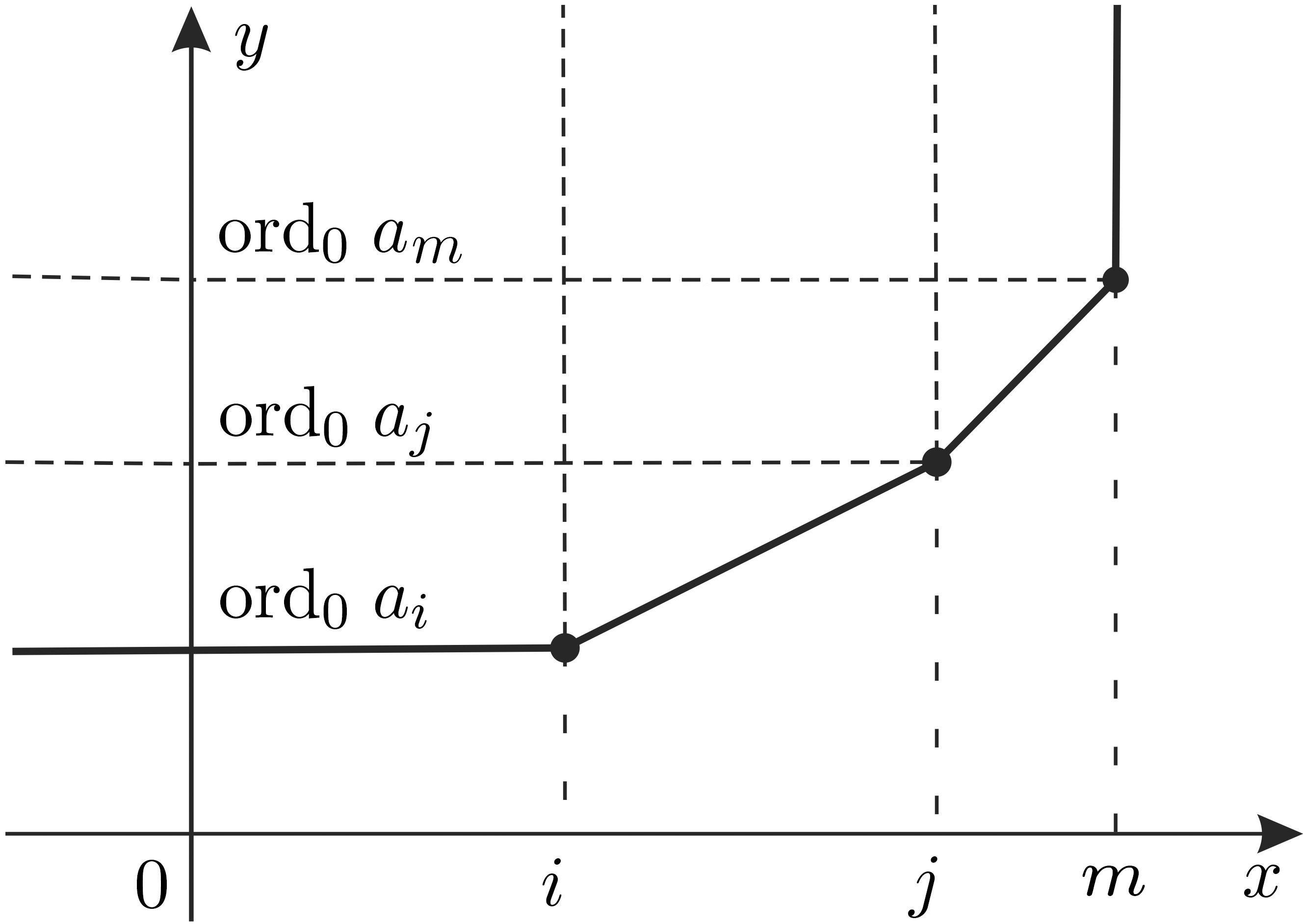}
\caption{\,The Newton polygon ${\cal N}(L)$ with two positive slopes $\displaystyle k_1=({\rm ord}_0\,a_j-{\rm ord}_0\,a_i)/(j-i),\;$ 
$\displaystyle k_2=({\rm ord}_0\,a_m-{\rm ord}_0\,a_j)/(m-j)$.}
\end{figure}

The result of Ramis has been further generalized by B.\,Malgrange and Y.\,Sibuya for a non-linear ODE of the general form (\ref{ADE}).
\medskip

{\bf Theorem 1} (Malgrange \cite{Mal}). {\it Let $\varphi\in{\mathbb C}[[z]]$ satisfy the equation $(\ref{ADE})$, that is 
$F(z,\Phi)=0$, where $\Phi=(\varphi,\delta\varphi,\ldots,\delta^m\varphi)$, and $\frac{\partial F}{\partial u_m}(z,\Phi)\ne0$. Then $\varphi$ is a power series of Gevrey order $1/k$, where $k$ is the least of all the positive slopes of the Newton polygon ${\cal N}(L_{\varphi})$ of a linear operator
$$
L_{\varphi}=\sum_{i=0}^m\frac{\partial F}{\partial u_i}(z,\Phi)\,\delta^i
$$ 
$($or $k=+\infty$, if ${\cal N}(L_{\varphi})$ has no positive slopes$)$.}
\medskip

The refinement of Theorem 1 belongs to Y.\,Sibuya \cite[App. 2]{Si}. This claims that $\varphi$ is a power series of the 
{\it exact} Gevrey order $1/k\in\{0,1/k_1,\ldots,1/k_s\}$, where $k_1<\ldots<k_s<\infty$ are all of the positive slopes 
of the Newton polygon ${\cal N}(L_{\varphi})$.
\medskip

In the paper we study {\it generalized} power series solutions of (\ref{ADE}) of the form 
\begin{eqnarray}\label{series}
\varphi=\sum_{n=0}^{\infty}c_nz^{s_n}, \qquad c_n\in{\mathbb C}, \qquad s_n\in{\mathbb C},
\end{eqnarray}
with the power exponents satisfying conditions
$$
0\leqslant{\rm Re}\,s_0\leqslant{\rm Re}\,s_1\leqslant\ldots, \qquad \lim_{n\rightarrow\infty}{\rm Re}\,s_n=+\infty
$$
(the latter, in particular, implies that a set of exponents having a fixed real part is finite).

Note that substituting the series (\ref{series}) into the equation (\ref{ADE}) makes sense, as only a finite number of terms in $\varphi$ contribute to any term of the form $cz^s$ in the expansion of $F(z,\Phi)=F(z,\varphi,\delta\varphi,\ldots,\delta^m\varphi)$ in powers of $z$. Indeed, $\delta^j\varphi=\sum_{n=0}^{\infty}c_ns_n^jz^{s_n}$ and an equation $s=s_{n_0}+s_{n_1}+\ldots+s_{n_l}$ has a finite number of solutions $(s_{n_0},s_{n_1},\ldots,s_{n_l})$, since $0\leqslant{\rm Re}\,s_n\rightarrow+\infty$. Furthermore, for any integer $N$ an inequality
${\rm Re}(s_{n_0}+s_{n_1}+\ldots+s_{n_l})\leqslant N$ has also a finite number of solutions, so that powers of $z$ in the expansion of $F(z,\Phi)$ can be ordered by the increasing of real parts. Thus, one may correctly define the notion of 
a formal solution of (\ref{ADE}) in the form of a generalized power series. In particular, the Painlev\'e III, V, VI equations are known to have such formal solutions (see \cite{Sh}, \cite{BG}, \cite{Gu}, \cite{Pa}).

For the generalized power series (\ref{series}) one may naturally define the {\it valuation}
$$
{\rm val}\,\varphi=s_0,
$$
and this is also well defined for any polynomial in $z,\varphi,\delta\varphi,\ldots,\delta^m\varphi$.

The main result of the paper is an analogue of the Maillet theorem (more precisely, of the Malgrange theorem) for 
generalized power series.
\medskip

{\bf Theorem 2.} {\it Let the generalized power series $(\ref{series})$ formally satisfy the equation $(\ref{ADE})$,
$\frac{\partial F}{\partial u_m}(z,\Phi)\ne0$, and for each $i=0,1,\ldots,m$ one have
\begin{eqnarray}\label{condition}
\frac{\partial F}{\partial u_i}(z,\Phi)=A_iz^{\lambda}+B_iz^{\lambda_i}+\ldots, \qquad 
{\rm Re}\,\lambda_i>{\rm Re}\,\lambda,
\end{eqnarray}
where not all $A_i$ equal zero. Let $k$ be the least of all the positive slopes of the Newton polygon ${\cal N}(L_{\varphi})$
$($or $k=+\infty$, if ${\cal N}(L_{\varphi})$ has no positive slopes$)$. Then for any sector $S$ of sufficiently small radius 
with the vertex at the origin and of the opening less than $2\pi$, the series 
$$
\sum_{n=0}^{\infty}\frac{c_n}{\Gamma(1+s_n/k)}\,z^{s_n}
$$
converges uniformly in $S$.}
\medskip

The Newton polygon of $L_{\varphi}$ in the case of the generalized power series $\varphi$ is defined similarly to 
the classical case, as the boundary curve of a convex hull of a union
$$
\bigcup_{i=0}^m\Bigl\{(x,y)\in{\mathbb R}^2\mid x\leqslant i,\;y\geqslant{\rm Re}\,{\rm val}\,\frac{\partial F}{\partial u_i}(z,\Phi)\Bigr\}.
$$

Note that $k=+\infty$ in Theorem 2 if and only if $A_m\ne0$. In this case $\varphi$ does converge in $S$, 
which has been already proved in \cite{GG} by the majorant method. Here we consider the case $k<+\infty$ using 
other known methods rather than the majorant one (of course, the convergence could also be proved by these methods). 

The first step (Sections 2, 3) consists of representing the generalized power series solution (\ref{series}) of (\ref{ADE}) 
by a multivariate Taylor series and follows from the ''grid-basedness'' of this formal solution. The latter means that 
there are $\nu_1,\ldots,\nu_r\in\mathbb C$ such that any term $cz^s$ of $\varphi$ is of the form 
$c(z^{\nu_1})^{m_1}\ldots(z^{\nu_r})^{m_r}=cz^{m_1\nu_1+\ldots+m_r\nu_r}$, for some 
$m_1,\ldots,m_r\in{\mathbb Z}_+$. An earlier result of this kind is due to Grigoriev, Singer \cite{GS} studying 
generalized power series solutions with real exponents. According to later results \cite{vdH}, \cite{ADH}, more 
general (real) transseries solutions of a polynomial ODE are actually also grid-based (see \cite{Ed} as an introduction 
to transseries). The second step (Section 4) consists of applying the implicit mapping theorem for Banach spaces of 
the obtained multivariate Taylor series.  

\section{Reduction of the ODE to a special form}

As we mentioned above, we prove Theorem 2 in the case $k<+\infty$. This implies $A_m=0$ in the expansion (\ref{condition}) of 
$\frac{\partial F}{\partial u_m}(z,\Phi)$. Let $0\leqslant p<m$ be such that $A_p\ne0$ and $A_i=0$ for all $i>p$. Then the minimal
positive slope $k$ of the Newton polygon ${\cal N}(L_{\varphi})$ is
\begin{eqnarray}\label{slope}
k=\min_{i>p}\frac{{\rm Re}\,\lambda_i-{\rm Re}\,\lambda}{i-p}.
\end{eqnarray}

First we will reduce the equation (\ref{ADE}) to a special form. This is provided by a transformation
$$
u=\sum_{n=0}^{\mu}c_nz^{s_n}+z^{s_{\mu}}v,
$$
with $\mu\geqslant0$ that will be chosen later. We use ideas of \cite{Mal} adapted to the case of generalized power series.  

The formal solution (\ref{series}) can be represented in the form
$$
\varphi=\sum_{n=0}^{\mu}c_nz^{s_n}+z^{s_{\mu}}\psi=\varphi_{\mu}+z^{s_{\mu}}\psi, \quad
{\rm val}\,\psi=s_{\mu+1}-s_{\mu}.
$$
Taking into consideration the equality $\delta(z^{s_{\mu}}\psi)=z^{s_{\mu}}(\delta+s_{\mu})\psi$, we have the relations
$$
\delta^i(z^{s_{\mu}}\psi)=z^{s_{\mu}}(\delta+s_{\mu})^i\psi, \qquad i=1,\ldots,m.
$$
Therefore, denoting 
$$
\Phi=(\varphi,\delta\varphi,\ldots,\delta^m\varphi)=\Phi_{\mu}+z^{s_{\mu}}\Psi, \qquad \Psi=(\psi_0,\psi_1,\ldots,\psi_m)=
(\psi,(\delta+s_{\mu})\psi,\ldots,(\delta+s_{\mu})^m\psi),
$$ 
and applying the Taylor formula to the relation $F(z,\Phi)=0$, we have
\begin{eqnarray}\label{Taylor}
0 & = & F(z,\Phi_{\mu}+z^{s_{\mu}}\Psi)=F(z,\Phi_{\mu})+z^{s_{\mu}}\sum_{i=0}^m\frac{\partial F}{\partial u_i}(z,\Phi_{\mu})\psi_i+\nonumber \\
  &  & +\frac12z^{2s_{\mu}}\sum_{i,j=0}^m\frac{\partial^2F}{\partial u_i\partial u_j}(z,\Phi_{\mu})\psi_i\psi_j+\ldots.
\end{eqnarray}

According to the assumptions of Theorem 2 and definition of the integer $p$, for each $i=0,1,\ldots,m$ the formal series
$\frac{\partial F}{\partial u_i}(z,\Phi)$ is of the form
$$
\frac{\partial F}{\partial u_i}(z,\Phi)=A_iz^{\lambda}+B_iz^{\lambda_i}+\ldots, \qquad {\rm Re}\,\lambda_i>{\rm Re}\,\lambda,
$$
where $A_i=0$ for all $i>p$ and $A_p\neq 0$. Define a polynomial
\begin{equation}
L(\xi)=A_0+A_1(\xi+s_{\mu})+\ldots+A_p(\xi+s_{\mu})^p\label{Euler}
\end{equation}
of degree $p$ and choose an integer $\mu\geqslant0$ such that the following three conditions i), ii), iii) hold.
\begin{itemize}
\item[i)] $L(\xi)\ne0\quad\forall\xi$ with ${\rm Re}\,\xi>0$;
\item[ii)] ${\rm Re}(s_{\mu+1}-s_{\mu})>0$;
\item[iii)] ${\rm Re}\,s_{\mu}>{\rm Re}\,\lambda+2(m-p)k$. 
\end{itemize}
Note that for each $i=0,1,\ldots,m$, the real part of
\begin{eqnarray*}
{\rm val}\left(\frac{\partial F}{\partial u_i}(z,\Phi)-\frac{\partial F}{\partial u_i}(z,\Phi_{\mu})\right)=
{\rm val}\biggl(z^{s_{\mu}}\sum_{j=0}^m\frac{\partial^2 F}{\partial u_i\partial u_j}(z,\Phi_{\mu})\psi_j+\ldots\biggr)
\end{eqnarray*}
is greater than ${\rm Re}\,s_{\mu}>{\rm Re}\,\lambda$. Therefore (see (\ref{condition}) for $\frac{\partial F}{\partial u_i}(z,\Phi)$),
\begin{eqnarray}\label{decomposition}
\frac{\partial F}{\partial u_i}(z,\Phi_{\mu})=A_iz^{\lambda}+\widetilde B_iz^{\tilde\lambda_i}+\ldots, \qquad {\rm Re}\,\tilde\lambda_i>{\rm Re}\,\lambda.
\end{eqnarray}
Moreover, for $i>p$ we have
\begin{eqnarray}\label{tildelambda}
{\rm Re}\,\tilde\lambda_i\geqslant{\rm Re}\,\lambda+(i-p)k.
\end{eqnarray}
Indeed, if $\widetilde B_iz^{\tilde\lambda_i}=B_iz^{\lambda_i}$, then the above inequality follows from the definition (\ref{slope}) of the slope $k$,
otherwise one has ${\rm Re}\,\tilde\lambda_i>{\rm Re}\,s_{\mu}$, and (\ref{tildelambda}) follows from the condition iii) above. Thus,
(\ref{decomposition}) and (\ref{tildelambda}) imply a decomposition
\begin{eqnarray}\label{LL'}
z^{-\lambda}\sum_{i=0}^m\frac{\partial F}{\partial u_i}(z,\Phi_{\mu})(\delta+s_{\mu})^i=L(\delta)+L'(z,\delta),
\end{eqnarray}
where the polynomial $L$ is defined by the formula (\ref{Euler}), and exponents $\alpha$ in the monomials $z^{\alpha}(\delta+s_{\mu})^i$ of
the operator $L'(z,\delta)$ satisfy the inequalities
$$
{\rm Re}\,\alpha>0, \qquad {\rm Re}\,\alpha\geqslant(i-p)k.
$$
From the relation (\ref{Taylor}), condition ii) and inequality ${\rm Re}\,s_{\mu}>{\rm Re}\,\lambda$ it follows that
\begin{eqnarray}\label{valF}
{\rm Re}\,{\rm val}\,F(z,\Phi_{\mu})>{\rm Re}(s_{\mu}+\lambda).
\end{eqnarray}
Finally, dividing the relation (\ref{Taylor}) by $z^{s_{\mu}+\lambda}$ and using (\ref{LL'}), (\ref{valF}), and the condition iii), 
we obtain the equality of the form
\begin{eqnarray}\label{psirelation}
L(\delta)\psi+L'(z,\delta)\psi+N(z,z^{\nu}\psi_0,z^{\nu}\psi_1,\ldots,z^{\nu}\psi_m)=0,
\qquad \nu=(m-p)k,
\end{eqnarray}
where the linear differential operators $L(\delta)$, $L'(z,\delta)$ are described above, and $N(z,u_0,u_1,\ldots,u_m)$ 
is a finite linear combination of monomials of the form
$$
z^{\beta}u_0^{q_0}u_1^{q_1}\ldots u_m^{q_m}, \qquad \beta\in{\mathbb C},\;{\rm Re}\,\beta>0,\;q_i\in{\mathbb Z}_+.
$$
Thus, the transformation $u=\varphi_{\mu}+z^{s_{\mu}}v$ reduces the equation (\ref{ADE}) to an equation
\begin{eqnarray}\label{reduced}
L(\delta)v+L'(z,\delta)v+N(z,z^{\nu}v,z^{\nu}(\delta+s_{\mu})v,\ldots,z^{\nu}(\delta+s_{\mu})^mv)=0,
\end{eqnarray}
with a formal solution $v=\psi$. 
\medskip

{\bf Remark 1.} The condition i) is not used for the above reduction, but we add it from the beginning as this will be used
in the sequel. In particular, i) implies that the coefficients $c_n$, $n\geqslant\mu+1$, of the formal solution $\psi=\sum_{n=\mu+1}^{\infty}c_nz^{s_n-s_{\mu}}$ 
of (\ref{reduced}) are uniquely determined by the coefficients $c_0,\ldots,c_{\mu}$ of the partial sum $\varphi_{\mu}$ of $\varphi$ 
(whereas some of $c_0,\ldots,c_{\mu}$ themselves may be free parameters).

\section{Representation of a generalized power series solution by a multivariate Taylor series}

Let us define an additive semi-group $G$ generated by a (finite) set consisting of the number $\nu$ and all the power exponents $\alpha$, $\beta$ of 
the variable $z$ containing in the monomials $z^{\alpha}(\delta+s_{\mu})^i$, $z^{\beta}u_0^{q_0}u_1^{q_1}\ldots u_m^{q_m}$ of $L'(z,\delta)$, 
$N(z,u_0,u_1,\ldots,u_m)$ respectively. Let $r_1,\ldots,r_l$ be generators of this semi-group, that is,
$$
G=\{m_1r_1+\ldots+m_lr_l \mid m_i\in{\mathbb Z}_+,\;\sum_{i=1}^lm_i>0\}, \qquad {\rm Re}\,r_i>0.
$$
As a consequence of the relation (\ref{psirelation}) and condition that the real parts of the power exponents $\alpha$, $\beta$ are positive, we have 
the following auxiliary lemma (details of the proof see in \cite[Lemma 2]{GG}).   
\medskip

{\bf Lemma 1.} {\it All the numbers $s_n-s_{\mu}$, $n\geqslant\mu+1$, belong to the additive semi-group $G$.}
\medskip

We may assume that the generators $r_1,\ldots,r_l$ of $G$ are linearly independent over $\mathbb Z$. This is provided by 
the following lemma \cite[Lemma 3]{GG}.
\medskip 

{\bf Lemma 2.} {\it There are complex numbers $\rho_1,\ldots,\rho_{\tau}$ linearly independent over $\mathbb Z$, such that all ${\rm Re}\,\rho_i>0$ and an additive 
semi-group $G'$ generated by them contains the above semi-group $G$ generated by $r_1,\ldots,r_l$.}
\medskip

In what follows, for the simplicity of exposition we assume that $G$ is generated by two numbers:
$$
G=\{m_1r_1+m_2r_2 \mid m_{1,2}\in{\mathbb Z}_+,\;m_1+m_2>0\}, \qquad {\rm Re}\,r_{1,2}>0.
$$
In the case of an arbitrary number $l$ of generators all constructions are analogous, only multivariate Taylor series in $l$ rather than in two variables 
are involved.

We should estimate the growth of the coefficients $c_n$ of the generalized power series
$$
\psi=\sum_{n=\mu+1}^{\infty}c_nz^{s_n-s_{\mu}},
$$
which satisfies the equality (\ref{psirelation}). According to Lemma 1, all the exponents $s_n-s_{\mu}$ belong to the semi-group $G$:
$$
s_n-s_{\mu}=m_1r_1+m_2r_2, \qquad (m_1,m_2)\in M\subseteq{\mathbb Z}_+^2\setminus\{0\},
$$
for some set $M$ such that the map $n\mapsto(m_1,m_2)$ is a bijection from ${\mathbb N}\setminus\{1,\ldots,\mu\}$ to $M$. Hence,
$$
\psi=\sum_{(m_1,m_2)\in M}c_{m_1,m_2}z^{m_1r_1+m_2r_2}=\sum_{(m_1,m_2)\in{\mathbb Z}_+^2\setminus\{0\}}c_{m_1,m_2}z^{m_1r_1+m_2r_2}
$$
(in the last series one puts $c_{m_1,m_2}=0$, if $(m_1,m_2)\not\in M$).

Now we define a natural linear map $\sigma:{\mathbb C}[[z^G]]\rightarrow{\mathbb C}[[z_1,z_2]]_*$ from 
the $\mathbb C$-algebra of generalized power series with exponents in $G$ to the $\mathbb C$-algebra of Taylor series in two 
variables without a constant term,
$$
\sigma: \sum_{\gamma=m_1r_1+m_2r_2\in G}a_{\gamma}z^{\gamma}\mapsto
\sum_{\gamma=m_1r_1+m_2r_2\in G}a_{\gamma}z_1^{m_1}z_2^{m_2}.
$$
As follows from the linear independence of the generators $r_1,r_2$ over $\mathbb Z$, 
$$
\sigma(\eta_1\eta_2)=\sigma(\eta_1)\sigma(\eta_2) \qquad\forall\eta_1,\eta_2\in{\mathbb C}[[z^G]],
$$
hence $\sigma$ is an isomorphism. The differentiation $\delta:{\mathbb C}[[z^G]]\rightarrow{\mathbb C}[[z^G]]$ naturally induces
a linear bijective map $\Delta$ of ${\mathbb C}[[z_1,z_2]]_*$ to itself,
$$
\Delta: \sum_{\gamma\in G}a_{\gamma}z_1^{m_1}z_2^{m_2}\mapsto\sum_{\gamma\in G}\gamma\,a_{\gamma}z_1^{m_1}z_2^{m_2},
$$
which clearly satisfies $\Delta\circ\sigma=\sigma\circ\delta$, so that the following commutative diagramme holds:
$$
\begin{array}{ccc}
{\mathbb C}[[z^G]] & \stackrel{\delta}{\longrightarrow} & {\mathbb C}[[z^G]] \\
\downarrow\lefteqn{\sigma} & & \downarrow\lefteqn{\sigma} \\
{\mathbb C}[[z_1,z_2]]_* & \stackrel{\Delta}{\longrightarrow} & {\mathbb C}[[z_1,z_2]]_*
\end{array}
$$

Thus we have the representation 
$$
\tilde\psi=\sigma(\psi)=\sum_{\gamma\in G}c_{\gamma}z_1^{m_1}z_2^{m_2}
$$ 
of the formal solution $\psi$ of (\ref{reduced}) by a multivariate Taylor series, where $c_{\gamma}=c_{m_1,m_2}$ for every $\gamma=m_1r_1+m_2r_2$. 
Now we apply the map $\sigma$ to the both sides of the equality (\ref{psirelation}), 
$$
L(\delta)\psi+L'(z,\delta)\psi+N(z,z^{\nu}\psi_0,z^{\nu}\psi_1,\ldots,z^{\nu}\psi_m)=0,
$$
and obtain a relation for $\tilde\psi$:
\begin{eqnarray}\label{psi0relation}
L(\Delta)\tilde\psi+\widetilde L(z_1,z_2,\Delta)\tilde\psi+\widetilde N(z_1,z_2,\,z_1^{k_1}z_2^{k_2}\tilde\psi,z_1^{k_1}z_2^{k_2}\tilde\psi_1,\ldots,
z_1^{k_1}z_2^{k_2}\tilde\psi_m)=0,
\end{eqnarray}
where $\tilde\psi_i=\sigma(\psi_i)=(\Delta+s_{\mu})^i\tilde\psi$, and
\begin{itemize}
\item[a)] $\widetilde L(z_1,z_2,\Delta)$ is a finite linear combination of monomials $z_1^{l_1}z_2^{l_2}(\Delta+s_{\mu})^i$ which satisfy 
$$
{\rm Re}(l_1r_1+l_2r_2)\geqslant(i-p)k,
$$
furthermore $\widetilde L(0,0,\Delta)\equiv0$;
\item[b)] $\widetilde N(z_1,z_2,u_0,\ldots,u_m)$ is a polynomial, and $\widetilde N(0,0,u_0,\ldots,u_m)\equiv0$;
\item[c)] $\nu=(m-p)k=k_1r_1+k_2r_2$.
\end{itemize}

\section{Banach spaces of multivariate Taylor series and the implicit mapping theorem}

In this section we conclude the proof of Theorem 2 establishing the corresponding required properties for the 
multivariate Taylor series 
$\tilde\psi$, which represents the generalized power series $\psi$ and satisfies the relation (\ref{psi0relation}).
We use the dilatation method based on the implicit mapping theorem for Banach spaces. This was originally used by
Malgrange for proving Theorem 1, as well as by C.\,Zhang \cite{Zh} in a further generalization of this theorem for 
$q$-difference-differential equations. 

Let us define the following Banach spaces $H^j$ of (formal) Taylor series in two variables without a constant term:
$$
H^j=\Bigl\{\eta=\sum_{\gamma\in G}a_{\gamma}z_1^{m_1}z_2^{m_2}\mid\sum_{\gamma\in G}\frac{|\gamma|^j}{|\Gamma(\gamma/k)|}\,|a_{\gamma}|<
+\infty\Bigr\}, \qquad j=0,1,\ldots,p,
$$ 
with the norm
$$
\|\eta\|_j=\sum_{\gamma\in G}\frac{|\gamma|^j}{|\Gamma(\gamma/k)|}\,|a_{\gamma}|.
$$
(The completeness of each $H^j$ is checked in a way similar to that how one checks the completeness of the space $l_2$;
see, for example, \cite[Ch. 6, \S4]{Die}.) One clearly has $H^p\subset H^{p-1}\ldots\subset H^0$ and
$$
\Delta+s_{\mu}: H^j\rightarrow H^{j-1}, \qquad j=1,\ldots,p.
$$
Therefore, $(\Delta+s_{\mu})^p$ maps $H^p$ to $H^0$, whereas $(\Delta+s_{\mu})^i$ for $i>p$ may map $H^p$ outside $H^0$. However,
linear operators $z_1^{l_1}z_2^{l_2}(\Delta+s_{\mu})^i$ with suitable $l_1$, $l_2$ possess the following property.
\medskip

{\bf Lemma 3.} {\it Let $l_1$, $l_2$ be such that ${\rm Re}(l_1r_1+l_2r_2)\geqslant(i-p)k$. Then  
$$
z_1^{l_1}z_2^{l_2}(\Delta+s_{\mu})^i: H^p\rightarrow H^0 
$$
is a continous linear operator.} 
\medskip

{\bf Proof.} Let $\eta=\sum_{\gamma\in G}a_{\gamma}z_1^{m_1}z_2^{m_2}\in H^p$. Then, by definition,
$$
\sum_{\gamma\in G}\frac{|\gamma+s_{\mu}|^p}{|\Gamma(\gamma/k)|}\,|a_{\gamma}|<+\infty,
$$
therefore
$$
a_{\gamma}=\frac{\Gamma(\gamma/k)}{(\gamma+s_{\mu})^p}\,A_{\gamma},
$$
where $\sum_{\gamma\in G}A_{\gamma}$ is an absolutely convergent series. Thus we have
\begin{eqnarray*}
z_1^{l_1}z_2^{l_2}(\Delta+s_{\mu})^i(\eta)&=&\sum_{\gamma\in G}(\gamma+s_{\mu})^i\,a_{\gamma}z_1^{m_1+l_1}z_2^{m_2+l_2}=
\sum_{\gamma>\gamma'}(\gamma-\gamma'+s_{\mu})^i\,a_{\gamma-\gamma'}\,z_1^{m_1}z_2^{m_2}=\\
&=&\sum_{\gamma>\gamma'}(\gamma-\gamma'+s_{\mu})^{i-p}\,\Gamma(\gamma/k-\gamma'/k)A_{\gamma-\gamma'}\,z_1^{m_1}z_2^{m_2},
\end{eqnarray*}
where we write $\gamma>\gamma'$ for $\gamma=m_1r_1+m_2r_2$, $\gamma'=l_1r_1+l_2r_2$, if $m_1\geqslant l_1$, $m_2\geqslant l_2$
and $(m_1,m_2)\ne(l_1,l_2)$.

To check the inclusion $z_1^{l_1}z_2^{l_2}(\Delta+s_{\mu})^i(\eta)\in H^0$, we should prove the absolute convergence of the series
\begin{eqnarray}\label{h0}
\sum_{\gamma>\gamma'}(\gamma-\gamma'+s_{\mu})^{i-p}\,\frac{\Gamma(\gamma/k-\gamma'/k)}{\Gamma(\gamma/k)}\,A_{\gamma-\gamma'}.
\end{eqnarray}
Since for $\gamma\rightarrow\infty$, ${\rm Re}\,\gamma>0$, and a fixed $\gamma'$, one has (see \cite[Ch. I, \S1.4]{BE}) 
$$
\frac{\Gamma(\gamma/k-\gamma'/k)}{\Gamma(\gamma/k)}\sim C(\gamma-\gamma')^{-\gamma'/k},
$$
a general term of the series (\ref{h0}) is equivalent to $C(\gamma-\gamma')^{i-p-\gamma'/k}A_{\gamma-\gamma'}$. This implies the
convergence of this series, as ${\rm Re}(i-p-\gamma'/k)\leqslant0$ under the assumptions of the lemma. We also have 
$$
\|z_1^{l_1}z_2^{l_2}(\Delta+s_{\mu})^i(\eta)\|_0\leqslant C_1\sum_{\gamma\in G}|A_{\gamma}|\leqslant C_2\|\eta\|_p,
$$
whence the continuity of $z_1^{l_1}z_2^{l_2}(\Delta+s_{\mu})^i$ follows. {\hfill $\Box$}
\medskip

{\bf Lemma 4.} {\it For any $\eta_1$, $\eta_2\in H^0$, one has $\eta_1\eta_2\in H^0$ and $\|\eta_1\eta_2\|_0\leqslant C\|\eta_1\|_0\|\eta_2\|_0$.}
\medskip

{\bf Proof.} Let $\eta_1=\sum_{\gamma\in G}a_{\gamma}z_1^{m_1}z_2^{m_2}$ and $\eta_2=\sum_{\gamma\in G}b_{\gamma}z_1^{m_1}z_2^{m_2}$. Then
$$
\eta_1\eta_2=\sum_{\gamma\in G}\Bigl(\sum_{\tilde\gamma<\gamma}a_{\tilde\gamma}b_{\gamma-\tilde\gamma}\Bigr)z_1^{m_1}z_2^{m_2},
$$ 
so that
$$
\|\eta_1\eta_2\|_0=\sum_{\gamma\in G}\frac{|\sum_{\tilde\gamma<\gamma}a_{\tilde\gamma}b_{\gamma-\tilde\gamma}|}{|\Gamma(\gamma/k)|}.
$$
Using the well known relation
$$
\frac{\Gamma(\tilde\gamma/k)\Gamma(\gamma/k-\tilde\gamma/k)}{\Gamma(\gamma/k)}=\int_0^1(1-t)^{\frac{\tilde\gamma}k-1}\,t^{\frac{\gamma-\tilde\gamma}k-1}\,dt,
$$
we have
$$
\left|\frac{\Gamma(\tilde\gamma/k)\Gamma(\gamma/k-\tilde\gamma/k)}{\Gamma(\gamma/k)}\right|\leqslant
\int_0^1(1-t)^{{\rm Re}\,\frac{\tilde\gamma}k-1}\,t^{{\rm Re}\,\frac{\gamma-\tilde\gamma}k-1}\,dt\leqslant C
$$
for any $\tilde\gamma<\gamma\in G$. Hence,
$$
\frac1{|\Gamma(\gamma/k)|}\leqslant\frac C{|\Gamma(\tilde\gamma/k)\Gamma(\gamma/k-\tilde\gamma/k)|},
$$
which implies
$$
\|\eta_1\eta_2\|_0\leqslant C\sum_{\gamma\in G}\sum_{\tilde\gamma<\gamma}\frac{|a_{\tilde\gamma}|}{|\Gamma(\tilde\gamma/k)|}\,
\frac{|b_{\gamma-\tilde\gamma}|}{|\Gamma(\gamma/k-\tilde\gamma/k)|}=C\|\eta_1\|_0\|\eta_2\|_0.
$$
{\hfill $\Box$}
\medskip

Now we conclude the proof of Theorem 2 by the implicit mapping theorem for Banach spaces which we recall below (see \cite[Th. 10.2.1]{Die}).
\medskip

{\it Let $\cal E$, $\cal F$, $\cal G$ be Banach spaces, $A$ an open subset of the direct product ${\cal E}\times{\cal F}$ and 
$h: A\rightarrow\cal G$ a continuously differentiable mapping. Consider a point $(x_0, y_0)\in A$ such that 
$h(x_0,y_0)=0$ and $\frac{\partial h}{\partial y}(x_0,y_0)$ is a bijective linear mapping from $\cal F$ to $\cal G$.

Then there are a neighbourhood $U_0\subset\cal E$ of the point $x_0$ and a unique continuous mapping
$g: U_0\rightarrow\cal F$ such that $g(x_0)=y_0$, $(x,g(x))\in A$, and $h(x,g(x))=0$ for any $x\in U_0$.}
\medskip

We will apply this theorem to the Banach spaces $\mathbb C$, $H^p$, $H^0$, and to the mapping $h:{\mathbb C}\times H^p\rightarrow H^0$ defined by
$$
h: (\lambda,\eta)\mapsto
L(\Delta)\eta+\widetilde L(\lambda z_1,\lambda z_2,\Delta)\eta+\widetilde N(\lambda z_1,\lambda z_2,\,\lambda^{k_1+k_2}z_1^{k_1}z_2^{k_2}\eta,\ldots,
\lambda^{k_1+k_2}z_1^{k_1}z_2^{k_2}(\Delta+s_{\mu})^m\eta),
$$
with $L$, $\widetilde L$, and $\widetilde N$ coming from (\ref{psi0relation}). This mapping is continuously differentiable by Lemmas 3, 4, moreover
$h(0,0)=0$ and $\frac{\partial h}{\partial\eta}(0,0)=L(\Delta)$ is a bijective linear mapping from $H^p$ to $H^0$. Therefore, there are $\rho>0$ and
$\eta_{\rho}\in H^p$ such that  
$$
L(\Delta)\eta_{\rho}+\widetilde L(\rho z_1,\rho z_2,\Delta)\eta_{\rho}+\widetilde N(\rho z_1,\rho z_2,\,\rho^{k_1+k_2}z_1^{k_1}z_2^{k_2}\eta_{\rho},\ldots,
\rho^{k_1+k_2}z_1^{k_1}z_2^{k_2}(\Delta+s_{\mu})^m\eta_{\rho})=0.
$$
Making the change of variables $(z_1,z_2)\mapsto(\frac{z_1}{\rho},\frac{z_2}{\rho})$, which induces an automorphism 
$\eta(z_1,z_2)\mapsto\eta(\frac{z_1}{\rho},\frac{z_2}{\rho})$ of ${\mathbb C}[[z_1,z_2]]_*$ commuting with $\Delta$, one can easily see that the 
above relation implies that the power series $\eta_{\rho}(\frac{z_1}{\rho},\frac{z_2}{\rho})$ satisfies the same equality 
(\ref{psi0relation}) as $\tilde\psi=\sum_{\gamma\in G}c_{\gamma}z_1^{m_1}z_2^{m_2}$ does. Hence, these two series coincide (the coefficients of a series satisfying (\ref{psi0relation}) are determined uniquely by this equality) and 
$$
\sum_{\gamma\in G}\frac{c_{\gamma}}{\Gamma(\gamma/k)}z_1^{m_1}z_2^{m_2}
$$
has a non-zero radius of convergence. This implies (substitute $z_1=z^{r_1}$, $z_2=z^{r_2}$ remembering that ${\rm Re}\,r_{1,2}>0$) the convergence of
the series
$$
\sum_{\gamma\in G}\frac{c_{\gamma}}{\Gamma(1+\gamma/k)}z^{\gamma}=\sum_{n=\mu+1}^{\infty}\frac{c_n}{\Gamma(1+(s_n-s_{\mu})/k)}z^{s_n-s_{\mu}}
$$
for any $z$ from a sector $S$ of sufficiently small radius with the vertex at the origin and of the opening less than $2\pi$, whence Theorem 2 follows.
\medskip

{\bf Remark 2.} From Theorem 2 one deduces the following estimate for the coefficients $c_n$ of the formal series solution (\ref{series}) of (\ref{ADE}): 
$$
|c_n|\leqslant A\,B^{{\rm Re}\,s_n}|\Gamma(1+s_n/k)|,
$$
for some $A,B>0$.

\end{document}